\documentclass[letterpaper,12pt,reqno]{amsart}
\usepackage[english]{babel}
\usepackage{amsthm,amsfonts,amssymb,amsmath}
\usepackage{graphicx}
\usepackage[noadjust]{cite}
\usepackage{fullpage}
\usepackage{verbatim}
\usepackage{todonotes}
\usepackage{ulem}

\newcommand{\C}{\mathbb{C}}

\newcommand{\Z}{\mathbb{Z}}
\newcommand{\es}[1]{\begin{equation}\begin{split}#1\end{split}\end{equation}}
\newcommand{\est}[1]{\begin{equation*}\begin{split}#1\end{split}\end{equation*}}

\newcommand{\tn}[1]{\textnormal{#1}}
\newtheorem*{theo*}{Theorem}
\newtheorem{theo}{Theorem}
\newtheorem{lemma}{Lemma}
\newtheorem{corol}{Corollary}

\newtheorem*{rem*}{Remark}
\newtheorem{conj}{Conjecture}
\newcommand{\pr}[1]{\left( #1\right)}

\newcommand{\sign}{\operatorname{sign}}

\begin{document}

\author{S. Al\.i Altu\u{g}, Sandro Bettin, Ian Petrow, Rishikesh, and Ian Whitehead}
\address{S. Al\.i Altu\u{g} -- Columbia University, Department of Mathematics, Room 517, 2990 Broadway, New York, NY 10027, USA}
\address{Sandro Bettin --  Centre de Recherches Math\'{e}matiques - Universit\'{e} de Montr\'{e}al, P.O. Box 6128, Centre-ville Station, Montr\'{e}al, QC, H3C 3J7, Canada}
\address{Ian Petrow -- Stanford University, Department of Mathematics, 450 Serra Mall, Building \#380, Stanford, CA, 94305, USA}
\address{Rishikesh -- Concordia University, 1455 de Maisonneuve Blvd. West, LB-901, Montr\'{e}al, QC, H3G 1M8, Canada}
\address{Ian Whitehead -- Columbia University, Department of Mathematics, Room 610, 2990 Broadway, New York, NY 10027, USA}

\title{A Recursion Formula for Moments of Derivatives of Random Matrix Polynomials}

\maketitle 

\begin{abstract} We give asymptotic formulae for random matrix
  averages of derivatives of characteristic polynomials over the
  groups $\tn{USp}(2N), \tn{SO}(2N)$ and $\tn{O}^-(2N)$. These
  averages are used to predict the asymptotic formulae for moments of
  derivatives of $L$-functions which arise in number theory. Each
  formula gives the leading constant of the asymptotic in terms of
  determinants of hypergeometric functions. We find a differential
  recurrence relation between these determinants which allows the
  rapid computation of the ($k+1$)-st constant in terms of the $k$-th
  and ($k-1$)-st. This recurrence is reminiscent of a Toda lattice
  equation arising in the theory of $\tau$-functions associated with
  Painlev\'{e} differential equations.     
\end{abstract}

\section{Introduction}

For over 50 years, mathematicians and physicists have used random
matrix theory to study a wide-ranging, growing list of probabilistic
phenomena. Particularly surprising are its applications in number
theory, where random matrices model the distribution of nontrivial
zeros of the Riemann zeta function. Random matrix theory now provides
far-reaching and widely believed conjectures for many questions in the
analytic theory of $L$-functions. 

Katz and Sarnak \cite{KS} give evidence that every family of $L$-functions falls into one of four symmetry types: unitary $\tn{U}(N)$, unitary symplectic $\tn{USp}(2N)$, even orthogonal $\tn{SO}(2N)$ and odd orthogonal $\tn{O}^-(2N)$. These symmetry types govern the
distribution of zeroes and special values in families. Using random
matrix models, Keating and Snaith \cite{KeatingSnaith} and Conrey,
Farmer, Keating, Rubinstein and Snaith \cite{CFKRS2} have produced
deep conjectures for estimating the integral moments of central
values in families of $L$-functions.    

The derivatives of $L$-functions are also of great interest, and are
the subject of this paper. A motivational example is Speiser's
theorem, which asserts that the Riemann hypothesis is equivalent to
the nonexistence of nonreal zeros of the derivative of the Riemann zeta function to the
left of the critical line, see e.g. \cite{SoundDuke}. Moreover, the
derivatives of $L$-functions control the order of vanishing at the
central point, which encodes important arithmetic and geometric
information. For example, according to the Birch and Swinnerton-Dyer
Conjecture, the order of vanishing of the $L$-function of an elliptic
curve over the rationals coincides with the arithmetic rank of the
curve. 

An $L$-function is modeled by the characteristic polynomial
$\Lambda_A$ of a random matrix $A$. Here we
compute \[M_k(\tn{G}(2N),m) := \int_{\tn{G}(2N)}
\pr{\Lambda^{(m)}_A(1)}^k\,dA,\] 
where $\tn{G}$ denotes $\tn{USp}$, $\tn{SO}$, or $\tn{O}^-$, and $dA$
is the Haar measure on $\tn{G}$. As $N\to\infty$, this models the
$k$th moment of $L^{(m)}(1/2)$ in a family of symmetry type $\tn{G}$. 

One can find the moments of $\Lambda^{(m)}_A(1)$ by differentiating
the corresponding shifted moment formulae, which are computed in
\cite{CFKRS2}. Conrey, Rubinstein and Snaith in \cite{CRS} develop a
faster method to compute the relevant averages in the unitary
case:\[\int_{\tn{U}(N)}{|\Lambda'_A(1)|}^{2k}\,d
A=b_kN^{k^2+2k}+O\pr{N^{k^2+2k-1}}.\]  
This leading constant $b_k$ is the same ``geometric constant''
appearing in conjectures for the asymptotic estimate of the $k$-th
moment of $\zeta'(s)$ in $t$ aspect. Conrey, Rubinstein and Snaith
describe $b_k$ in terms of a $k \times k$ determinant of $I$-Bessel
functions, and are able to compute $b_k$ numerically for $k\leq 15$.

Forrester and Witte in \cite{FW, FW2002} find a surprising expression for
these determinants of $I$-Bessel functions in terms of solutions to
Painlev\'{e} III$^\prime$ differential equations.  Explicitly, the formula in \cite{CRS} is as follows:
 \[b_k = (-1)^k \sum_{h = 0 }^k \binom{k}{h}
 \left.\left(\frac{d}{dt}\right)^{h+k} \left(e^{-t} t^{-k^2/2} \det_{k
       \times k} \left( I_{k+i-j}(2 \sqrt{t})\right) \right)
 \right|_{t=0} \]  
where $I_\ell(x)$ is the modified Bessel function of the first kind. One then defines
\[ \tau_k(t) := 2^{-k(k-1)} t^{-k^2/2} \det_{k \times k} \left( I_{k+i-j}(2 \sqrt{t})\right).\]
Forrester and Witte~\cite{FW,FW2002} find that this $\tau_k(t)$
(denoted $\tau[k](t)$ in \cite[section 4]{FW2002}) is in fact the
Okamoto $\tau$ function associated with the Painlev\'{e} III$^\prime$
differential equation: \[ (ty'')^2 + y'(4y'-1)(y-ty')-\frac{1}{4}k^2 =
0.\] This nonlinear second order differential equation has a solution
with certain boundary data (see \cite{FW2002}) given in terms of
$\tau_k(t)$ by the formula \[ y= \sigma_{\text{III},k}(t) = - t
\frac{d}{dt} \log \left( e^{-t/4} t^{k^2}
  \tau_k\left(\frac{t}{4}\right)\right).\]  
Specifying boundary conditions, one can quickly compute
$\sigma_{\text{III},k}(t)$ from the differential equation and recover
$\tau_k(t)$ via the equation 
\[
\tau_k(t)=\exp\left(-\int_0^{4t}(\sigma_{\text{III},k}(s)+k^2-\frac{s}{4})\frac{ds}{s}\right).\]
This expression allows a much faster computation of the constants
$b_k$.  

The goal of the present paper is to extend these results to the other symmetry types relevant to $L$-functions: $\tn{USp}(2N)$, $\tn{SO}(2N)$ and $\tn{O}^-(2N)$. 

Employing similar techniques to \cite{CRS}, we obtain analogous results, where the role of the $I$-Bessel functions above is here played by hypergeometric functions,

\est{
g_m(u)&=\frac{1}{2\pi i}\oint_{|w|=1}\frac{e^{w+\frac{u}{w^2}}}{w^{m+1}} \, dw\\
&=\frac{1}{\Gamma(m+1)}{}_{0}F_{2}\pr{;\frac{m}2+1,\frac{m+1}{2};\frac{u}{4}},\\
}
for $u\in\C$ and $m\in\Z$. For negative $m,$ interpret the above expression as the limit. The role of the $\tau$-function is played
by
\begin{equation}
  \label{taufunction}
\mathcal{T}_{k,\ell}(u):=\det_{k\times k}\pr{g_{2i-j+\ell}(u)},  
\end{equation}
for $k\geq0$, $\ell\in\Z$ and $u\in\C$, where, here and in the following, the indices $i$ and $j$ of the matrix in the determinant range from $1$ to $k$. In the context of Theorems \ref{MTSY}, \ref{MTSO} and \ref{MTOO} below, the $\ell$ appearing here takes values $0,-1,0$ respectively. We now state theorems for the symplectic, special orthogonal, and negative orthogonal cases: 

\begin{theo}\label{MTSY}
We have
\begin{equation}\label{meqsy}
M_k(\tn{USp}(2N),2) = b_k(\tn{USp}(2N),2)\cdot(2N)^{\frac{k^2+5k}{2}} + O(N^{\frac{k^2+3k}{2}})
\end{equation}
where 
\est{
b_k(\tn{USp}(2N),2)= 2^{-\frac{k^2+5k}{2}} \frac{d^k}{du^{k}}\left.\pr{e^u \mathcal{T}_{k,0}(2u)}\right|_{u=0}.
}
\end{theo}

\begin{theo}\label{MTSO}
We have
\begin{equation}\label{meqso}
M_k(\tn{SO}(2N),2) = b_k(\tn{SO}(2N),2)\cdot(2N)^{\frac{k^2+3k}{2}} + O(N^{\frac{k^2+k}{2}})
\end{equation}
where 
\est{
b_k(\tn{SO}(2N),2)= 2^{-\frac{k^2+k}{2}} \frac{d^k}{du^{k}}\left.\pr{e^u \mathcal{T}_{k,-1}(2u)}\right|_{u=0}.
}
\end{theo}

\begin{theo}\label{MTOO}
We have
\begin{equation}\label{meqoo}
M_k(\tn{O}^-(2N),3) = b_k(\tn{O}^-(2N),3)\cdot (2N)^{\frac{k^2+5k}{2}} + O(N^{\frac{k^2+3k}{2}})
\end{equation}
where 
\est{
b_k(\tn{O}^-(2N),3)= 3\cdot 2^{-\frac{k^2+3k}{2}} \frac{d^k}{du^{k}}\left.\pr{e^u \mathcal{T}_{k,0}(2u)}\right|_{u=0}.
}
\end{theo}
Note that we find above that $b_k(\tn{O}^-(2N),2)= 3\cdot 2^k \cdot b_k(\tn{USp}(2N),2)$.

We are naturally led to consider the second and third derivatives of characteristic polynomials in the above theorems instead of the first derivative due to root number considerations. Indeed, if $A$ is a unitary matrix, the characteristic polynomial satisfies the functional equation
 \[\Lambda_A(s) = (-s)^N (\det A)^{-1} \overline{\Lambda}_{A}(s^{-1}),\]
 where $\overline{f}(s)=\overline{f(\overline{s})}$. When $A$ is in
 $\tn{USp}(2N)$ or $\tn{SO}(2N)$ then $\det A$ is constantly equal to
 $+1$ and one has a simple expression for $\Lambda'_A(1)$ in terms of $\Lambda_A(1)$, which can be used to compute the moments of the derivative via partial integration. Thus the moments of $\Lambda''_A(1)$ give the next novel information. When $A \in \tn{O}^-(2N)$ we have $\det A= -1$, and thus $\Lambda_A(1) = 0$. In this case $\Lambda'_A(1)$ plays the role which $\Lambda_A(1)$ plays in the other families, $\Lambda''_A(1)$ has a simple expression in terms of $\Lambda'_A(1)$, and therefore one considers moments of $\Lambda'''_A(1)$. The same reasoning carries over to the families of $L$-functions having each of the aforementioned symmetry types. 

The method of proof of Theorems \ref{MTSY}, \ref{MTSO} and \ref{MTOO} can be generalized to higher-order derivatives easily. It suffices to expand the binomial in Lemma~\ref{lder}, and otherwise proceed as in the given proofs of Theorems \ref{MTSY}, \ref{MTSO} and \ref{MTOO}. 

We can use Theorems~\ref{MTSY},~\ref{MTSO} and~\ref{MTOO} to give conjectures for moments of derivatives of $L$-functions at $s=1/2$ with the above symmetry types. For example, the quadratic Dirichlet $L$-functions ordered by conductor form a symplectic family and thus we make the following conjecture:

\begin{conj}
Let $\mathcal{D}(X):=\{ |d|<X,d \tn{ fundamental discriminant}\}$, and $L(s,\chi_d)$ denote the quadratic Dirichlet L-function of fundamental discriminant $d$.
The average value of the second derivative of quadratic Dirichlet L-functions at the central point is
\est{
\frac{1}{|\mathcal{D}(X)|}\sum_{d\in \mathcal{D}(X)}L''(1/2, \chi_d)^k\sim a_k\cdot b_k(\tn{USp}(2N),2)\cdot (\log X)^{\frac{k^2+5k}{2}}
}
\end{conj}
Here $a_k$ is a well understood arithmetic constant depending on the functional equations in the family:
\[
a_k=\prod_{p\text{ prime}}\frac{(1-1/p)^{k(k+1)/2}}{1+1/p}\pr{
\frac{(1-1/\sqrt{p})^{-k}+(1+1/\sqrt{p})^{-k}}{2}+\frac1p}
\]

(see, for example~\cite[1.3.5]{CFKRS}). This is the same arithmetic constant appearing in moment conjectures for $L(1/2, \chi_d)$ without the derivative. It is not predicted by random matrix calculations.

Our Theorems \ref{MTSY}, \ref{MTSO} and \ref{MTOO} along with the results of \cite{CRS} allow similar conjectures to be made for any family of $L$-functions. 

Under closer examination, the determinants $\mathcal{T}_{k,\ell}(u)$
defined by \eqref{taufunction} and appearing in Theorems \ref{MTSY}, \ref{MTSO} and \ref{MTOO} exhibit a surprisingly rich structure. Our Theorem \ref{recurrence} is a differential recurrence relation which allows much faster computation of the constants $b_k(\tn{G}(N),m)$:

\begin{theo}\label{recurrence}Let $k\in\mathbb{Z}_{>0}$, $\ell\in\mathbb{Z}$. Then
\begin{equation}\label{lem1}
\mathcal{T}_{k+1,\ell}(u)\mathcal{T}_{k-1,\ell}(u)=2\left(u\mathcal{T}_{k,\ell}(u)\mathcal{T}_{k,\ell}''(u)+\mathcal{T}_{k,\ell}(u)\mathcal{T}_{k,\ell}'(u)-u\left(\mathcal{T}_{k,\ell}'(u)\right)^2\right).
\end{equation}
\end{theo}

This recurrence relation closely resembles a Toda lattice equation for the Okamoto $\tau$-function associated with a Painlev\'{e} differential equation, see \cite[Theorem 2]{Okamoto87}.  Such Toda lattice equations are at the heart of the $\tau$-function theory of Painlev\'{e} equations, and are used by Forrester and Witte \cite{FW,FW2002} to connect determinants the of $I$-Bessel functions found by \cite{CRS} to the Painlev\'{e} III${}^\prime$ equation.  It would be very interesting to determine whether or not there exists a differential equation arising from our formula \eqref{lem1} which plays the role for symplectic and orthogonal types that Painlev\'{e} III${}^\prime$ plays for unitary symmetry.

Ultimately, one hopes to obtain formulae for the complex moments of characteristic polynomials in the various symmetry types.  In the case of the undifferentiated moment conjectures, it has been found that the geometric constants $g_k$ can be expressed in a simple form in terms of Barnes $G$-functions, which are well-defined for complex values of $k$, see~\cite{CF}.  A project for the future would be to see if there exists a similar expression for the geometric coefficients $b_k$ studied in this paper for moments of derivatives of $L$-functions.

By computing the expressions in Theorems \ref{MTSY}, \ref{MTSO} and \ref{MTOO} directly, one can obtain the values of $b_k$ up to $k\approx 10$. By using Theorem \ref{recurrence} we do much better: running SAGE for about an hour on a machine with 4 gigabytes of RAM we computed the first 200 values of $b_k(\tn{USp}(2N),2)$. In section \ref{numbers} we give a table with the first 10 values of $b_k$ for each symmetry type. To give an example, we have 
\[ b_{10}({\tn{USp}(2N)})= \frac { 
47 \cdot
1553 \cdot
1787 \cdot
73709 \cdot
152825093 
}
{
2 ^{ 62 } \cdot
3 ^{ 34 } \cdot
5 ^{ 17 } \cdot
7 ^{ 10 } \cdot
11 ^{ 5 } \cdot
13 ^{ 5 } \cdot
17 ^{ 4 } \cdot
19 ^{ 3 } \cdot
23 ^{ 2 } \cdot
29 \cdot
31 \cdot
37 
},\] \[ b_{10}(\tn{SO}(2N))=\frac { 
25171 \cdot
7695491 \cdot
57668937071891 
}
{
2 ^{ 45 } \cdot
3 ^{ 29 } \cdot
5 ^{ 15 } \cdot
7 ^{ 9 } \cdot
11 ^{ 5 } \cdot
13 ^{ 5 } \cdot
17 ^{ 3 } \cdot
19 ^{ 3 } \cdot
23 ^{ 2 } \cdot
29 \cdot
31 \cdot
37 
},\] \[b_{10}({\tn{O}^-(2N)})= \frac { 
47 \cdot
1553 \cdot
1787 \cdot
73709 \cdot
152825093 
}
{
2 ^{ 52 } \cdot
3 ^{ 33 } \cdot
5 ^{ 17 } \cdot
7 ^{ 10 } \cdot
11 ^{ 5 } \cdot
13 ^{ 5 } \cdot
17 ^{ 4 } \cdot
19 ^{ 3 } \cdot
23 ^{ 2 } \cdot
29 \cdot
31 \cdot
37 
}.\]

\textbf{Acknowledgements}: 
We would like to thank the organizers and staff of the American Mathematical Society's Mathematics Research Community program on Arithmetic Statistics in Snowbird, Utah, who made our work possible. We thank the staff of the American Institute of Mathematics for their hospitality. Particular thanks to Brian Conrey for his guidance throughout this project.

\section{Some lemmas}

We recall the definitions of the relevant spaces of matrices: $\tn{USp}(2N)$ is the subgroup of $2N \times 2N$ unitary matrices $M$ with $M^t\left(\begin{array}{cc} 0 & I \\ -I & 0 \end{array}\right)M=\left(\begin{array}{cc} 0 & I \\ -I & 0 \end{array}\right)$, where I denotes the $N\times N$ identity matrix. $\tn{SO}(2N)$ and $\tn{O}^-(2N)$ are the subsets ($\tn{O}^-$ is not a subgroup) of orthogonal matrices with determinant $1$ and $-1$, respectively. These compact spaces admit Haar measures which we normalize so that the volume of each space is $1$. A matrix in $\tn{USp}(2N)$ or $\tn{SO}(2N)$ has characteristic polynomial of the form $\Lambda(x)=\prod_{n=1}^N (1-e^{i\theta_n}x)(1-e^{-i\theta_n}x)$, and a matrix in $\tn{O}^-(2N)$ has characteristic polynomial of the form $\Lambda(x)=(1-x)(1+x)\prod_{n=1}^{N-1} (1-e^{i\theta_n}x)(1-e^{-i\theta_n}x)$, with $\theta_n \in \mathbb{R}$.

The shifted moments of $\Lambda$ are defined as follows:
\est{
I(\tn{G}(2N);z_1,\dots,z_k):=\int_{\tn{G}(2N)}\Lambda(z_1)\cdots \Lambda(z_k)\,dA.
}
Conrey, Farmer, Keating, Rubinstein, and Snaith \cite{CFKRS2} use Weyl integration formula (see, for example, Theorem 8.60 of \cite{Kn}) to compute the following shifted moment formulae, which are the starting point for our work (note that we have corrected a typo in \cite[ 4.9]{CFKRS2}):

\begin{lemma}
\cite[3.36]{CFKRS2} Assume that $\alpha_1,\dots,\alpha_k$ are complex
numbers with $|\alpha_i|<1$ for $i=1,\dots k$. Then 
\est{
I(&\tn{USp}(2N);e^{-\alpha_1},\dots,e^{-\alpha_k})=\\
&=\frac{(-1)^{k(k-1)/2}2^k}{(2\pi i)^kk!}\oint_{|w_1|=1}\!\!\!\!\cdots\oint_{|w_k|=1}\frac{\prod_{1\leq i <j \leq k}(w_i^2-w_j^2)^2\prod_{j=1}^kw_j}{\prod_{1\leq i,j\leq k}(w_j^2-\alpha_i^2)} \times\\
&\hspace{2.7cm}\times e^{N\sum_{j=1}^k(w_j-\alpha_j)}\prod_{1\leq m\leq \ell\leq k}(1-e^{-w_m-w_\ell})^{-1}\,dw_1\cdots dw_k.\\
}
\end{lemma}

\begin{lemma}
\cite[4.43]{CFKRS2} Assume that $\alpha_1,\dots,\alpha_k$ are complex
numbers with $|\alpha_i|<1$ for $i=1,\dots k$. Then 
\est{
I(&\tn{SO}(2N);e^{-\alpha_1},\dots,e^{-\alpha_k})=\\
&=\frac{(-1)^{k(k-1)/2}2^k}{(2\pi i)^kk!}\oint_{|w_1|=1}\!\!\!\!\cdots\oint_{|w_k|=1}\frac{\prod_{1\leq i <j \leq k}(w_i^2-w_j^2)^2\prod_{j=1}^kw_j}{\prod_{1\leq i,j\leq k}(w_j^2-\alpha_i^2)} \times\\
&\hspace{2.7cm}\times e^{N\sum_{j=1}^k(w_j+\alpha_j)}\prod_{1\leq m< \ell< k}(1-e^{-w_m-w_\ell})^{-1}\,dw_1\cdots dw_k.\\
}
\end{lemma}

\begin{lemma}
\cite[ 4.9]{CFKRS2} Assume that $\alpha_1,\dots,\alpha_k$ are complex numbers with $|\alpha_i|<1$ for $i=1,\dots k$. Then
\est{
I(&\tn{O}^-(2N);e^{-\alpha_1},\dots,e^{-\alpha_k})=\\
&=\frac{(-1)^{k(k-1)/2}2^k}{(2\pi i)^kk!}\oint_{|w_1|=1}\!\!\!\!\cdots\oint_{|w_k|=1}\frac{\prod_{1\leq i <j \leq k}(w_i^2-w_j^2)^2\prod_{j=1}^k\alpha_j}{\prod_{1\leq i,j\leq k}(w_j^2-\alpha_i^2)} \times\\
&\hspace{2.7cm}\times e^{N\sum_{j=1}^k(w_j+\alpha_j)}\prod_{1\leq m\leq \ell\leq k}(1-e^{-w_m-w_\ell})^{-1}\,dw_1\cdots dw_k.\\
}
\end{lemma}

In fact, we will use approximate versions of the lemmas above, which follow immediately from the fact that $(1-e^{-x})^{-1}=x^{-1}+O(1)$. To simplify notation, we denote Vandermonde determinants as follows:
\est{
\Delta(w):=\det_{k\times k}(w_i^{j-1})=\prod_{1\leq i<j\leq k}\pr{w_i-w_j}
}
and we write $w^2=(w_i^2)_{1\leq i\leq k}$ for any $w=(w_i)_{1\leq i\leq k}\in\C^k$.
Then we have:

\begin{corol}\label{mcsy}
Assume that $\alpha_1,\dots,\alpha_k$ are complex numbers such that $|\alpha_j|\ll\frac 1{N}$ for $j=1\dots k$. Then
\est{
& I(\tn{USp}(2N);e^{-\alpha_1},\dots,e^{-\alpha_k})=\\ & =\frac{(-1)^{k(k-1)/2}}{(2\pi i)^kk!} \left( \oint \cdots \oint \frac{\Delta(w)\Delta(w^2)e^{N\sum_j(w_j-\alpha_j)}}{\prod_{i,j}(w_j^2-\alpha_i^2)}dw_1\dots dw_k \right)(1+O(N^{-1})).\\
}
\end{corol}

\begin{corol}\label{mcso}
Assume that $\alpha_1,\dots,\alpha_k$ are complex numbers such that $|\alpha_j|\ll\frac 1{N}$ for $j=1\dots k$. Then
\est
{
& I(\tn{SO}(2N);e^{-\alpha_1},\dots,e^{-\alpha_k})=\\ & =\frac{(-1)^{k(k-1)/2}2^k}{(2\pi i)^kk!}\left( \oint \cdots\oint \frac{\Delta(w)\Delta(w^2)(\prod_j{w_j}) e^{N\sum_j(w_j+\alpha_j)}}{\prod_{i,j}(w_j^2-\alpha_i^2)}\,dw_1\dots dw_k \right) (1+O(N^{-1})).
}
\end{corol}

\begin{corol}\label{mcoo}
Assume that $\alpha_1,\dots,\alpha_k$ are complex numbers such that
$|\alpha_j|\ll\frac 1{N}$ for $j=1\dots k$. Then 
\est
{
& I(\tn{O}^-(2N);e^{-\alpha_1},\dots,e^{-\alpha_k})=\\ & =\frac{(-1)^{k(k-1)/2}2^k}{(2\pi i)^kk!}\left( \oint \cdots \oint \frac{\Delta(w)\Delta(w^2)(\prod_j{\alpha_j})e^{N\sum_j (w_j+\alpha_j)}}{\prod_{i,j}(w_j^2-\alpha_i^2)}\,dw_1\dots dw_k \right) (1+O(N^{-1})).\\
}
\end{corol}

The following formula will appear in the proofs of Theorems~\ref{MTSY} and~\ref{MTSO}, below. We use it only when $m=2$, but the general statement is a starting point for computing moments of the $m$th derivative.

\begin{lemma}\label{lder}
For $m\geq0$, we have
\est{
&\frac{d^m}{d\alpha_1^m}\cdots\frac{d^m}{d\alpha_k^m}\left.\frac{e^{-N\sum_{i=1}^k\alpha_i}}{\prod_{1\leq i,j\leq k}(w_j^2-\alpha_i^2)}\right|_{\alpha_1=\dots=\alpha_k=0}=\\
&\hspace{4cm}=\Bigg(\sum_{\ell=0}^m{{m}\choose{\ell}}(-N)^{m-\ell}\sum_{\substack{i_1+\cdots +i_k=\ell,\\ i_j\ \text{even}}}\prod_{j=1}^k\frac{i_j!}{w_j^{i_j+2}}\Bigg)^k\\
}
\end{lemma}
\begin{proof}
We have
\est{
\frac{d^m}{d\alpha^m}\frac{e^{-N\alpha}}{\prod_{1\leq j\leq k}(w_j^2-\alpha^2)}&=\sum_{\ell=0}^m{{m}\choose{\ell}}(-N)^{m-\ell}e^{-N\alpha}\sum_{i_1+\cdots i_k=m}\prod_{j=1}^k\frac{d^{i_j}}{d\alpha^{i_j}}\frac{1}{(w_j^2-\alpha^2)}\\
}
and
\est{
\left.\frac{d^q}{d\alpha^q}\left(\frac{1}{w^2-\alpha^2}\right)\right|_{\alpha=0}
=\begin{cases}
0 & \text{if }q \text{ is odd}\\ 
q!/w^{q+2} & \text{if }q \text{ is even.}
\end{cases}
}
Therefore,
\est{
\left. \frac{d^m}{d\alpha^m}\frac{e^{-N\alpha}}{\prod_{1\leq j\leq k}(w_j^2-\alpha^2)} \right|_{\alpha=0} &=\sum_{\ell=0}^m{{m}\choose{\ell}}(-N)^{m-\ell}e^{-N\alpha}\sum_{\substack{i_1+\cdots +i_k=\ell,\\ i_j\ \text{even}}}\prod_{j=1}^k{\frac{i_j!}{w_j^{i_j+2}}}\\
}
and the lemma follows.
\end{proof}

Our proofs will also rely upon Vandermonde determinants of differential operators:
\est{\Delta\pr{\frac{d}{dx}}:=\prod_{1\leq i<j\leq k}\pr{\frac{d}{dx_i}-\frac{d}{dx_j}}=\det_{ k\times k}\pr{\frac{d^{j-1}}{dx_i^{j-1}}}.
}
We give two lemmas on computing with these--Lemma~\ref{fld} is a direct consequence of the definition, but we prove Lemma~\ref{ufl} in detail. 
\begin{lemma}\label{fld}
Let $f_1(x),\dots,f_k(x)$ be $k-1$ times differentiable. Then
\est{
\Delta\pr{\frac{d}{dx}}\prod_{i=1}^kf_i(x_i)=\det_{k\times k}\pr{f_i^{(j-1)}(x_i)}.
}
\end{lemma}
\begin{lemma}\label{ufl}
Let $f(x,y)$ be $k-1$ times differentiable in $x$ and $y$. Then 
\est{
\Delta\pr{\frac{d}{dx}}\Delta\pr{\frac{d}{dy}}\left.\prod_{i=1}^k f(x_i,y_i)\right|_{\substack{x_1=\ldots =x_k=X,\\y_1=\ldots=y_k=Y}}=k!\det_{k\times k}\pr{\frac{d^{i+j-2}}{dX^{i-1}dY^{j-1}}f(X,Y)}
}
\end{lemma}
\begin{proof}
By Lemma~\ref{fld} we have
\est{
\Delta\pr{\frac{d}{dx}}\prod_{i=1}^kf(x_i,y_i)=\det_{k\times k}\pr{\frac{d^{j-1}}{dx_i^{j-1}}f(x_i,y_i)}=\sum_{\mu}\sign(\mu)\prod_{i=1}^k\frac{d^{\mu(i)-1}}{dx_i^{\mu(i)-1}}f(x_i,y_i)
}
where the sum runs over the permutations $\mu\in S_k$. Applying Lemma~\ref{fld} again, we find
\est{
\Delta\pr{\frac{d}{dy}}\Delta\pr{\frac{d}{dx}}\prod_{i=1}^kf(x_i,y_i)&=\sum_{\mu}\sign(\mu)\Delta\pr{\frac{d}{dy}}\prod_{i=1}^k\frac{d^{\mu(i)-1}}{dx_i^{\mu(i)-1}}f(x_i,y_i)\\
&=\sum_{\mu}\sign(\mu)\det_{k\times k}\pr{\frac{d^{\mu(i)+j-2}}{dx_i^{\mu(i)-1}dy_i^{j-1}}f(x_i,y_i)}\\
}
and so 
\est{
\Delta\pr{\frac{d}{dy}}\Delta\pr{\frac{d}{dx}}\left.\prod_{i=1}^kf(x_i,y_i)\right|_{\substack{x_i=X,\\y_i=Y}}&=\sum_{\mu}\sign(\mu)\det_{k\times k}\pr{\frac{d^{\mu(i)+j-2}}{dX^{\mu(i)-1}dY^{j-1}}f(X,Y)}.
}
Now, we may rearrange the rows of the matrix \[\pr{\frac{d^{\mu(i)+j-2}}{dX^{\mu(i)-1}dY^{j-1}}f(X,Y)}\] to obtain \[\pr{\frac{d^{i+j-2}}{dX^{i-1}dY^{j-1}}f(X,Y)};\] doing so cancels out the $\sign(\mu)$ attached to the determinant, and we reach the desired formula.
\end{proof}

Our final lemma is a recursion for determinants discovered by Lewis Carroll. It will be used in the proof of Theorem~\ref{recurrence}.

\begin{lemma}\label{LCL}
Let $A$ be an $k\times k$ matrix, and let $A{{a_1,\dots,a_r}\choose{b_1,\dots,b_s}}$ denote the matrix $A$ with rows $a_1,\dots,a_r$ and the columns $b_1,\dots,b_s$ removed. Then
\est{
\det A{{i}\choose{i}} \cdot \det A{{j}\choose{j}}-\det A{{i}\choose{j}} \cdot \det A{{j}\choose{i}}=\det A \cdot \det A{{i,j}\choose{i,j}}.
}
\end{lemma}

\section{Proof of the theorems}
A simple calculation shows that
\est{
\frac{d^m}{d\alpha_1^m}\cdots\frac{d^m}{d\alpha_k^m}&\left.I(\tn{G}(2N);e^{-\alpha_1},\dots,e^{-\alpha_k})\right|_{\alpha_1=\dots\alpha_k=0}\\
&=(-1)^{mk}\int_{\tn{G}(2N)}\pr{\sum_{j=0}^m \left \{ \begin{matrix} m \\ j \end{matrix} \right \} \Lambda^{(j)}(1)}^k\,dA,
}
where $\left \{ \begin{matrix} m \\ j \end{matrix} \right \}$ denotes a Stirling number of the second kind. It follows that
\est{
M_{k}(\tn{G}(2N),m)=\frac{d^m}{d\alpha_1^m}\cdots\frac{d^m}{d\alpha_k^m}&\left.I(\tn{G}(2N);e^{-\alpha_1},\dots,e^{-\alpha_k})\right|_{\alpha_1=\cdots=\alpha_k=0}(1+O(N^{-1})).
}
Thus we may proceed towards Theorems~\ref{MTSY},~\ref{MTSO}, and~\ref{MTOO} by differentiating the formulae found in Corollaries~\ref{mcsy},~\ref{mcso}, and~\ref{mcoo}, respectively. The asymptotics remain valid after differentiating because they are uniform in $\alpha_i$.

\begin{proof}[Proof of Theorem~\ref{MTSY}]
From Corollary~\ref{mcsy} and the above argument, we know that the $k$th moment of $\Lambda''$ is asymptotically 
\est{
M_{k}(\tn{USp}(2N),2)=\frac{(-1)^{k(k-1)/2}}{k!} \tilde{M}_k(\tn{USp}(2N),2)(1+O(N^{-1}))
}
where
\begin{multline*}
 \tilde{M}_k(\tn{USp}(2N),2) = \\
\frac{d^{2k}}{d\alpha_1^2\cdots d\alpha_k^2}\left. \frac{1}{(2\pi
    i)^k}\oint \cdots \oint
  \frac{\Delta(w)\Delta(w^2)e^{N\sum_j(w_j-\alpha_j)}}{\prod_{i,j}(w_j^2-\alpha_i^2)}dw_1\dots
  dw_k \right|_{\alpha_1\dots=\alpha_k=0}.
\end{multline*}

We apply Lemma~\ref{lder} with $m=2$, finding
\es{\label{dbb}
\frac{d^{2k}}{d\alpha_1^2\cdots
  d\alpha_k^2}\left.\frac{e^{-N\sum_{i=1}^k\alpha_i}}{\prod_{1\leq
      i,j\leq
      k}(w_j^2-\alpha_i^2)}\right|_{\alpha_1=\dots=\alpha_k=0}&=\pr{\prod_{j=1}^k\frac{1}{w_j^{2k}}}
\Bigg(N^2+2\sum_{j=1}^k\frac{1}{w_j^2}\Bigg)^k\\ 
&=\frac{d^k}{dt^{k}}\left.\pr{\prod_{j=1}^k\frac{1}{w_j^{2k}}}\exp\Bigg(tN^2+2t\sum_{j=1}^k\frac{1}{w_j^2}\Bigg)\right|_{t=0}.
}
This allows us to write 
\est{
& \tilde{M}_k(\tn{USp}(2N),2) = \frac{d^k}{dt^{k}}\left. \frac{e^{tN^2}}{(2\pi i)^k} \oint \cdots \oint \Delta(w)\Delta(w^2)\exp\pr{\sum_{j=1}^k Nw_j+\frac{2t}{w_j^2}}\,\frac{d w_1}{w_1^{2k}}\cdots \frac{dw_k}{w_k^{2k}}\right|_{t=0}.
}
We now replace the Vandermonde determinants in this expression with Vandermonde determinants of differential operators. Observe that 
\est{
\Delta(w^2)=\Delta\pr{\frac d{dL}}\left.e^{\sum_{i=1}^k w_i^2L_i}\right|_{L_i=0}.
}
and also
\est{
\Delta(w)\cdot e^{\sum_{i=1}^k Nw_i}= \Delta\pr{\frac d{dM}}\left.e^{\sum_{i=1}^kw_iM_i}\right|_{M_i=N}
}
This implies that we can compute the integral in the above formula as 
\est{
& \frac{1}{(2\pi i)^k} \oint \cdots \oint \Delta(w)\Delta(w^2)\exp\pr{\sum_{j=1}^k Nw_j+\frac{2t}{w_j^2}}\,\frac{d w_1}{w_1^{2k}}\cdots \frac{dw_k}{w_k^{2k}} \\
& = \Delta\pr{\frac d{dL}}\Delta\pr{\frac d{dM}} \frac{1}{(2\pi i)^k} \oint \cdots \oint \exp\pr{\sum_{j=1}^k L_j w_j^2+M_jw_j+\frac{2t}{w_j^2}}\,\frac{d w_1}{w_1^{2k}}\cdots \frac{dw_k}{w_k^{2k}}\bigg|_{\substack{L_j=0,\\M_j=N}} \\
&=\Delta\pr{\frac d{dL}}\Delta\pr{\frac d{dM}}\prod_{j=1}^k \pr{\frac{1}{2\pi i}\oint_{|w|=1}\exp\pr{L_j w^2+M_j w+\frac{2t}{w^2}}\frac{d w}{w^{2k}}}\bigg|_{\substack{L_j=0,\\M_j=N}} \\
&= k! \det_{k \times k}\pr{\frac{d^{i+j-2}}{dL^{i-1}dM^{j-1}} \, \frac{1}{2\pi i}\oint_{|w|=1}\exp\pr{L w^2+M w+\frac{2t}{w^2}}\frac{d w}{w^{2k}} \Bigg|_{\substack{L=0,\\M=N}} }
}
where the last equality is by Lemma~\ref{ufl}. Using the fact that 
\est{
& \frac{d^{i+j-2}}{dL^{i-1}dM^{j-1}} \, \frac{1}{2\pi i}\oint_{|w|=1}\exp\pr{L w^2+M w+\frac{2t}{w^2}}\frac{d w}{w^{2k}} \Bigg|_{\substack{L=0,\\M=N}} \\ 
& = \frac{1}{2\pi i} \oint \frac{\exp(Nw+\frac{2t}{w^2})}{w^{2k-2i-j+3}}dw \\
& = \frac{N^{2k-2i-j+2}}{2\pi i} \oint \frac{\exp(w+\frac{2tN^2}{w^2})}{w^{2k-2i-j+3}}dw, 
}
we obtain a simplified formula for $\tilde{M}_k(\tn{USp}(2N),2)$:
\est{
\tilde{M}_k(\tn{USp}(2N),2) 
&= k! N^{k(k+1)/2} \frac{d^k}{dt^{k}} \left. \left( e^{tN^2} \det_{k\times k} \pr{\frac{1}{2\pi i} \oint \frac{\exp(w+\frac{2tN^2}{w^2})}{w^{2k-2i-j+3}}dw}\right) \right|_{t=0} \\
&=(-1)^{k(k-1)/2} k! N^{k(k+1)/2} \frac{d^k}{dt^{k}} \left. \left( e^{tN^2} \det_{k\times k} \pr{\frac{1}{2\pi i} \oint \frac{\exp(w+\frac{2tN^2}{w^2})}{w^{2i-j+1}}dw}\right) \right|_{t=0} \\
&=(-1)^{k(k-1)/2} k! N^{k(k+5)/2} \frac{d^k}{du^{k}} \left. \left( e^u \det_{k\times k} \pr{\frac{1}{2\pi i} \oint \frac{\exp(w+\frac{2u}{w^2})}{w^{2i-j+1}}dw}\right) \right|_{u=0}
}
where we have interchanged columns of the matrix to obtain the second line and set $u=tN^2$ to obtain the third. The theorem follows.
\end{proof}

\begin{proof}[Proof of Theorem~\ref{MTSO}]
We can proceed in the same way as in the proof of Theorem~\ref{MTSY}, starting from Corollary~\ref{mcso}.
\end{proof}

\begin{proof}[Proof of Theorem~\ref{MTOO}]
This time we consider the third derivative. From~\eqref{dbb} we have
\est{
\frac{d^{3k}}{d\alpha_1^3\cdots d\alpha_k^3}\left.\frac{\prod_{i=1}^k\alpha_i e^{N\sum_{i=1}^k\alpha_i}}{\prod_{1\leq i,j\leq k}(w_j^2-\alpha_i^2)}\right|_{\alpha_1=\dots=\alpha_k=0}&=3\frac{d^{2k}}{d\alpha_1^2\cdots d\alpha_k^2}\left.\frac{e^{N\sum_{i=1}^k\alpha_i}}{\prod_{1\leq i,j\leq k}(w_j^2-\alpha_i^2)}\right|_{\alpha_1=\dots=\alpha_k=0}\\
&=3\prod_{j=1}^k\frac{1}{w_j^{2k}}\Bigg(N^2+2\sum_{j=1}^k\frac{1}{w_j^2}\Bigg)^k\\
&=3\frac{d^k}{dt^{k}}\left.\prod_{j=1}^k\frac{1}{w_j^{2k}}\exp\Bigg(tN^2+2t\sum_{j=1}^k\frac{1}{w_j^2}\Bigg)\right|_{t=0}
}
and the theorem follows from Corollary~\ref{mcoo} in the same way.
\end{proof}

\begin{proof}[Proof of Theorem~\ref{recurrence}]

We begin by proving a two-variable version of the recurrence relation. Let
\est{
\tilde g_m(x,y)&:=\frac{1}{2\pi i}\int_{|z|=1}\frac{e^{xz+\frac{y}{z^2}}}{z^{m+1}} \, dz,\\
\tilde{\mathcal {T}}_{k,\ell}(x,y)&:=\det_{k\times k}\pr{\tilde g_{2i-j+\ell}(x,y)},\\
}
for $x,y\in\C$, $k\geq0$ and $\ell\in\Z$. Using Lemma \ref{LCL}, we will show that:
\begin{equation}\label{partial recurrence}
\tilde{\mathcal {T}}_{k+1,\ell}(x,y)\tilde{\mathcal {T}}_{k-1,\ell}(x,y)=\tilde{\mathcal {T}}_{k,\ell}(x,y)\frac{\partial^2}{\partial x\partial y}\tilde{\mathcal {T}}_{k,\ell}(x,y)-\frac{\partial}{\partial x}\tilde{\mathcal {T}}_{k,\ell}(x,y)\frac{\partial}{\partial y}\tilde{\mathcal {T}}_{k,\ell}(x,y).
\end{equation}
Let $A_{k, \ell}$ denote the matrix of $\tilde{\mathcal {T}}_{k,\ell}(x,y)$, i.e. 
\est{
A_{k, \ell}=\left(\begin{array}{cccc}
\tilde g_{1+\ell} & \tilde g_{\ell} & \cdots & \tilde g_{2-k+\ell} \\
\tilde g_{3+\ell} & \tilde g_{2+\ell} & \cdots & \tilde g_{4-k+\ell} \\
\cdots & \cdots &  & \cdots \\
\tilde g_{2k-1+\ell} & \tilde g_{2k-2+\ell} & \cdots & \tilde g_{k+\ell} 
\end{array}\right).
}
Observe that 
\est{
\frac{\partial }{\partial x}\tilde g_{m}(x,y)&=\tilde g_{m-1}(x,y,2t),\\
\frac{\partial }{\partial y}\tilde g_{m}(x,y)&=\tilde g_{m+2}(x,y,2t).
}
We now compute the partial derivatives of $\tilde{\mathcal {T}}_{k,\ell}(x,y)$. Expanding the derivative by columns, we obtain:
\est{
\frac{\partial}{\partial x} \tilde{\mathcal {T}}_{k,\ell} =  & \left|\begin{array}{cccc}
{\partial \tilde{g}_{1+\ell}}/{\partial x} & \tilde g_{\ell} & \cdots & \tilde g_{2-k+\ell} \\
{\partial \tilde{g}_{3+\ell}}/{\partial x} & \tilde g_{2+\ell} & \cdots & \tilde g_{4-k+\ell} \\
\cdots & \cdots &  & \cdots \\
{\partial \tilde{g}_{2k-1+\ell}}/{\partial x} & \tilde g_{2k-2+\ell} & \cdots & \tilde g_{k+\ell} 
\end{array}\right| + 
\left|\begin{array}{cccc}
\tilde g_{1+\ell} & {\partial \tilde{g}_{\ell}}/{\partial x} & \cdots & \tilde g_{2-k+\ell} \\
\tilde g_{3+\ell} & {\partial \tilde{g}_{2+\ell}}/{\partial x} & \cdots & \tilde g_{4-k+\ell} \\
\cdots & \cdots &  & \cdots \\
\tilde g_{2k-1+\ell} & {\partial \tilde{g}_{2k-1+\ell}}/{\partial x} & \cdots & \tilde g_{k+\ell} 
\end{array}\right| + \\
& \cdots +
\left|\begin{array}{cccc}
\tilde g_{1+\ell} & \tilde g_{\ell} & \cdots & {\partial \tilde{g}_{2-k+\ell}}/{\partial x} \\
\tilde g_{3+\ell} & \tilde g_{2+\ell} & \cdots & {\partial \tilde{g}_{4-k+\ell}}/{\partial x} \\
\cdots & \cdots &  & \cdots \\
\tilde g_{2k-1+\ell} & \tilde g_{2k-2+\ell} & \cdots & {\partial \tilde{g}_{k+\ell}}/{\partial x}
\end{array}\right|.
}
All terms but the last one in this sum vanish. Thus $\frac{\partial}{\partial x} \tilde{\mathcal {T}}_{k,\ell}=\det A_{k+1, \ell}{{k+1}\choose{k}}$, where $A{{a_1,\dots,a_r}\choose{b_1,\dots,b_s}}$ denotes the matrix $A$ with rows $a_1,\dots,a_r$ and the columns $b_1,\dots,b_s$ removed, as in Lemma \ref{LCL}. 

Similarly, expanding the partial derivative by rows, we have $\frac{\partial}{\partial y} \tilde{\mathcal {T}}_{k,\ell}=\det A_{k+1, \ell}{{k}\choose{k+1}}$.  Finally, $\frac{\partial^2}{\partial x \partial y} \tilde{\mathcal {T}}_{k,\ell}=\det A_{k+1, \ell}{{k}\choose{k}}$. Thus Equation \ref{partial recurrence} is an immediate consequence of Lemma \ref{LCL}.

To return to our original functions $g_m(u)$ and $\mathcal{T}_{k, \ell}(u)$, we observe that a simple change of variables gives $\tilde g(x,y)=x^mg_m(x^2y)$
for any $x\neq0$. Further, removing a factor of $x^{2i+\ell}$ from each row $i$ and $x^{-j}$ from each column $j$ in the determinant, we have 
\begin{equation}\label{lem6}
\tilde {\mathcal {T}}_{k,\ell}(x,y)=x^{\frac{k(k+1)}{2}+k\ell} {\mathcal {T}}_{k,\ell}(x^2y).
\end{equation}
If we set $u=x^2y$, then by the chain rule we have $\frac{\partial }{\partial x}=2xy\frac{d}{du}$, $\frac{\partial }{\partial y}=x^2\frac{d}{du}$,
and so
\est{
& x^{\frac{k(k+1)}{2}+kl}\mathcal{T}_{k,\ell}(u)\frac{\partial^2}{\partial x\partial y}\pr{x^{\frac{k(k+1)}{2}+kl}\mathcal{T}_{k,\ell}(u)} \\
&=x^{k(k+1)+2k\ell+1}\pr{\left(\frac{k(k+1)}{2}+k\ell+2\right)\mathcal{T}_{k,\ell}(u)\mathcal{T}_{k,\ell}'(u)+2u\mathcal{T}_{k,\ell}(u)\mathcal{T}_{k,\ell}''(u)}
}
and
\est{
& \frac{\partial}{\partial x}\pr{x^{\frac{k(k+1)}{2}+k\ell}\mathcal{T}_{k,\ell}(u)}\frac{\partial }{\partial y}\pr{{x^{\frac{k(k+1)}{2}+k\ell}}\mathcal{T}_{k,\ell}(u)} \\
&=x^{k(k+1)+2k\ell+1}\pr{\pr{\frac{k(k+1)}{2}+k\ell}\mathcal{T}_{k,\ell}(u)\mathcal{T}_{k,\ell}'(u)+2u \mathcal{T}'_{k,\ell}(u)^2}.
}
Thus the theorem follows from~\eqref{lem6} and~\eqref{partial recurrence}.

\end{proof}

\section{Numerical values}\label{numbers}
Below are the first several values for the constant $b_k(\tn{USp}(2N,2))$. See Theorem \ref{MTSY}. 
\begin{equation*} b_{ 1 } =\frac { 1 
}
{
2 \cdot
3 
}
 \end{equation*}
\begin{equation*} b_{ 2 } =\frac { 
19 
}
{
2 ^{ 4 } \cdot
3 ^{ 2 } \cdot
5 \cdot
7 
}
 \end{equation*}
\begin{equation*} b_{ 3 } =\frac { 
487 
}
{
2 ^{ 7 } \cdot
3 ^{ 5 } \cdot
5 ^{ 2 } \cdot
7 \cdot
11 
}
 \end{equation*}
\begin{equation*} b_{ 4 } =\frac { 
59 \cdot
197 
}
{
2 ^{ 13 } \cdot
3 ^{ 8 } \cdot
5 ^{ 2 } \cdot
7 ^{ 2 } \cdot
11 \cdot
13 
}
 \end{equation*}
\begin{equation*} b_{ 5 } =\frac { 
174290791 
}
{
2 ^{ 19 } \cdot
3 ^{ 10 } \cdot
5 ^{ 5 } \cdot
7 ^{ 3 } \cdot
11 ^{ 2 } \cdot
13 \cdot
17 \cdot
19 
}
 \end{equation*}
\begin{equation*} b_{ 6 } =\frac { 
3373 \cdot
1670407 
}
{
2 ^{ 25 } \cdot
3 ^{ 14 } \cdot
5 ^{ 6 } \cdot
7 ^{ 3 } \cdot
11 ^{ 3 } \cdot
13 ^{ 2 } \cdot
17 \cdot
19 \cdot
23 
}
 \end{equation*}
\begin{equation*} b_{ 7 } =\frac { 
37 \cdot
83 \cdot
2203 \cdot
3571457 
}
{
2 ^{ 32 } \cdot
3 ^{ 19 } \cdot
5 ^{ 9 } \cdot
7 ^{ 6 } \cdot
11 ^{ 3 } \cdot
13 ^{ 3 } \cdot
17 ^{ 2 } \cdot
19 \cdot
23 
}
 \end{equation*}
\begin{equation*} b_{ 8 } =\frac { 
61 \cdot
595351 \cdot
11423948521 
}
{
2 ^{ 42 } \cdot
3 ^{ 23 } \cdot
5 ^{ 11 } \cdot
7 ^{ 7 } \cdot
11 ^{ 4 } \cdot
13 ^{ 4 } \cdot
17 ^{ 2 } \cdot
19 ^{ 2 } \cdot
23 \cdot
29 \cdot
31 
}
 \end{equation*}
\begin{equation*} b_{ 9 } =\frac { 
53 \cdot
16646765854629827113 
}
{
2 ^{ 53 } \cdot
3 ^{ 29 } \cdot
5 ^{ 13 } \cdot
7 ^{ 9 } \cdot
11 ^{ 5 } \cdot
13 ^{ 4 } \cdot
17 ^{ 3 } \cdot
19 ^{ 3 } \cdot
23 ^{ 2 } \cdot
29 \cdot
31 
}
 \end{equation*}
\begin{equation*} b_{ 10 } =\frac { 
47 \cdot
1553 \cdot
1787 \cdot
73709 \cdot
152825093 
}
{
2 ^{ 62 } \cdot
3 ^{ 34 } \cdot
5 ^{ 17 } \cdot
7 ^{ 10 } \cdot
11 ^{ 5 } \cdot
13 ^{ 5 } \cdot
17 ^{ 4 } \cdot
19 ^{ 3 } \cdot
23 ^{ 2 } \cdot
29 \cdot
31 \cdot
37 
}
 \end{equation*}

 \vspace{0.2in}
 
Below are several values for the constant $b_k(\tn{SO}(2N),2)$. See Theorem \ref{MTSO}.

\begin{equation*} b_{ 1 } =1 \end{equation*}
\begin{equation*} b_{ 2 } =\frac { 
7 
}
{
2 \cdot
3 \cdot
5 
}
 \end{equation*}
\begin{equation*} b_{ 3 } =\frac { 
2 \cdot
13 
}
{
3 ^{ 4 } \cdot
5 \cdot
7 
}
 \end{equation*}
\begin{equation*} b_{ 4 } =\frac { 
17 \cdot
5987 
}
{
2 ^{ 6 } \cdot
3 ^{ 5 } \cdot
5 ^{ 2 } \cdot
7 ^{ 2 } \cdot
11 \cdot
13 
}
 \end{equation*}
\begin{equation*} b_{ 5 } =\frac { 
157 \cdot
17519 
}
{
2 ^{ 9 } \cdot
3 ^{ 8 } \cdot
5 ^{ 4 } \cdot
7 ^{ 2 } \cdot
11 \cdot
13 \cdot
17 
}
 \end{equation*}
\begin{equation*} b_{ 6 } =\frac { 
22273664659 
}
{
2 ^{ 15 } \cdot
3 ^{ 12 } \cdot
5 ^{ 5 } \cdot
7 ^{ 3 } \cdot
11 ^{ 2 } \cdot
13 ^{ 2 } \cdot
17 \cdot
19 
}
 \end{equation*}
\begin{equation*} b_{ 7 } =\frac { 
116228886131 
}
{
2 ^{ 14 } \cdot
3 ^{ 14 } \cdot
5 ^{ 8 } \cdot
7 ^{ 5 } \cdot
11 ^{ 3 } \cdot
13 ^{ 2 } \cdot
17 \cdot
19 \cdot
23 
}
 \end{equation*}
\begin{equation*} b_{ 8 } =\frac { 
36774351481263481 
}
{
2 ^{ 28 } \cdot
3 ^{ 19 } \cdot
5 ^{ 9 } \cdot
7 ^{ 6 } \cdot
11 ^{ 4 } \cdot
13 ^{ 3 } \cdot
17 \cdot
19 ^{ 2 } \cdot
23 \cdot
29 
}
 \end{equation*}
\begin{equation*} b_{ 9 } =\frac { 
71 \cdot
103 \cdot
223 \cdot
661 \cdot
1069 \cdot
134437 
}
{
2 ^{ 34 } \cdot
3 ^{ 25 } \cdot
5 ^{ 11 } \cdot
7 ^{ 6 } \cdot
11 ^{ 4 } \cdot
13 ^{ 4 } \cdot
17 ^{ 3 } \cdot
19 ^{ 2 } \cdot
23 \cdot
29 
}
 \end{equation*}
\begin{equation*} b_{ 10 } =\frac { 
25171 \cdot
7695491 \cdot
57668937071891 
}
{
2 ^{ 45 } \cdot
3 ^{ 29 } \cdot
5 ^{ 15 } \cdot
7 ^{ 9 } \cdot
11 ^{ 5 } \cdot
13 ^{ 5 } \cdot
17 ^{ 3 } \cdot
19 ^{ 3 } \cdot
23 ^{ 2 } \cdot
29 \cdot
31 \cdot
37 
}
 \end{equation*}

\vspace{0.2in}

We omit a table of values for the odd orthogonal case.  Recall that by Theorems \ref{MTSY} and \ref{MTOO}, $b_k(\tn{O}^-(2N),3)= 3\cdot 2^k \cdot b_k(\tn{USp}(2N),2)$ so that these values are given in terms of the above table for the symplectic case.

\vspace{12pt}

\end{document}